# A New Pivot Selection Algorithm for Symmetric Indefinite Factorization Arising in Quadratic Programming with Block Constraint Matrices

**Duangpen Jetpipattanapong\* and Gun Srijuntongsiri**
School of Information, Computer, and Communication Technology, Sirindhorn International Institute of Technology, Thammasat University, Thailand
\*Author for correspondence; e-mail: duangpenj@gmail.com

**ABSTRACT**

Quadratic programming is a class of constrained optimization problem with quadratic objective functions and linear constraints. It has applications in many areas and is also used to solve nonlinear optimization problems. This article focuses on the equality constrained quadratic programs whose constraint matrices are block diagonal. Using the direct solution method, we propose a new pivot selection algorithm for the factorization of the Karush-Kuhn-Tucker (KKT) matrix for this problem that maintains the sparsity and stability of the problem. Our experiments show that our pivot selection algorithm appears to produce no fill-ins in the factorization of such matrices. In addition, we compare our method with MA57 and find that the factors produced by our algorithm are sparser.

## 1. INTRODUCTION

Quadratic programming is a class of nonlinear programming. It has applications in many areas such as financial, signal processing, modeling, and simulation [1-8]. It is also the subproblem of sequential quadratic programming for solving nonlinear programs [9]. There are many methods to solve quadratic programs; they can be classified as either direct or iterative methods. Our work focuses on direct methods. These methods solve quadratic programs by directly solving the resulting linear systems using matrix factorization. Iterative methods, on the other hand, have two well-known classes of methods: active set and interior point methods. Unlike direct methods, the running time of iterative methods depend on the required accuracy of the solutions. Direct methods for solving quadratic programming typically involves symmetric indefinite factorization of the Karush-Kuhn-Tucker (KKT) matrix [10]. Symmetric indefinite factorization generally requires pivoting to maintain stability. Many pivot selection techniques have been proposed such as Bunch-Parlett [11], Bunch-Kaufman [12], and bounded Bunch-Kaufman method (BBK) [13]. Moreover, when the KKT matrix is sparse, the choice of pivots also affects the sparsity of the resulting factors, which in turn affects the time needed to solve the linear system. Hence, choosing suitable pivots that both maintain stability and preserve sparsity is not trivial. Duff et al. [14] use relative pivot tolerance for stability and apply the minimum degree algorithm [15] to consider a 2-by-2 pivot for the sparsity for sparse symmetric indefinite factorization. Olaf and Klaus [16] propose Supernode-Bunch-Kaufman pivoting method, which applies the Bunch-Kaufman pivot selection algorithm for the sparse case, supplemented by pivot perturbation techniques. Duff et al. introduce the multifrontal method for sparse symmetric indefinite factorization [17, 18]. The multifrontal approach is widely used in many sparse direct solver, for example, MA57 and MUMPS [19, 20].

In this article, we investigate equality-constrained quadratic programs whose constraint matrices are block diagonal. Such problems often arise in practice when different groups of variables are independent but variables in the same group must satisfy some constraints. Using a direct method, we propose a pivot selection algorithm for this type of quadratic programs. By exploiting the known structure of the quadratic program, our algorithm can efficiently identify the pivot candidates that can maintain the sparsity of the factors. We use the condition number of each pivot candidates as part of the information for pivot selection in order to also maintain stability. Our experiments show that our algorithm maintains sparsity and stability better than MA57. In Section 2, we define the block constraint quadratic programming problem, briefly explain symmetric indefinite factorization and describe our algorithm. Section 3 shows our experiment and the results. Finally, the conclusion is in Section 4.

## 2. MATERIALS AND METHODS
### 2.1 Quadratic Programs with Block Diagonal Constraint Matrices

This article considers equality constraint quadratic programs whose constraint matrices are block diagonal. The problem is as follows:
$$\min_{x \in \mathbb{R}^n} \frac{1}{2} x^T H x + c^T x, \quad (1)$$
subject to $Ax = e$

where $x \in \mathbb{R}^n, H \in \mathbb{R}^{n \times n}, c \in \mathbb{R}^n, A \in \mathbb{R}^{m \times n}, e \in \mathbb{R}^m, m < n$, and the constraint matrix $A$ is of the form
$$A = \begin{bmatrix} A_1 & 0 & \cdots & 0 \\ 0 & A_2 & \cdots & 0 \\ \vdots & \vdots & \ddots & \vdots \\ 0 & 0 & \cdots & A_N \end{bmatrix},$$
where $A_i \in \mathbb{R}^{m_i \times n_i}$ ($i = 1, 2, ..., N$), $m_i < n_i$, and $N$ is the number of diagonal blocks in $A$. Note that $\sum_{i=1}^N m_i = m$ and $\sum_{i=1}^N n_i = n$. Assume that $A$ has full row rank. Recall from the first-order necessary conditions that, for $x^*$ to be a solution of (1), there must be $x^*$ and $\lambda^*$ satisfying
$$\begin{bmatrix} H & -A^T \\ A & 0 \end{bmatrix} \begin{bmatrix} x^* \\ \lambda^* \end{bmatrix} = \begin{bmatrix} -c \\ e \end{bmatrix} \quad (2)$$
[10]. The above system of equations can be rewritten to a more useful form of Karush-Kuhn-Tucker (KKT) system
$$\begin{bmatrix} H & A^T \\ A & 0 \end{bmatrix} \begin{bmatrix} -p \\ \lambda^* \end{bmatrix} = \begin{bmatrix} g \\ h \end{bmatrix}, \quad (3)$$
where $p = x^* - x, g = c + Hx$, and $h = Ax - e$. The matrix in (3) is known as the KKT matrix.

### 2.2 Symmetric Indefinite Factorization

To solve the KKT system in (3), note that since the KKT matrix is symmetric indefinite, we cannot use Cholesky factorization to factorize it. Instead, we can perform symmetric indefinite factorization [21]. Let $K$ be the KKT matrix, a symmetric indefinite factorization of $K$ is in the following form
$$P^T K P = L B L^T, \quad (4)$$
where $L$ is a unit lower triangular matrix, $B$ is a block diagonal matrix with block dimension equal to 1 or 2, and $P$ is a permutation matrix. The permutation matrix $P$ is chosen to maintain numerical stability of the computation. In case $K$ is large and sparse, $P$ is chosen to also maintain the sparsity in $L$ in addition to maintaining the stability. After factorization, back and forward substitutions are used to compute the solution of (3) by the following steps:

(i) Solve $z : Lz = P^T \begin{bmatrix} g \\ h \end{bmatrix}$.

(ii) Solve $\hat{z} : B\hat{z} = z$.

(iii) Solve $\bar{z} : L^T \bar{z} = \hat{z}$.

(iv) Set: $\begin{bmatrix} -p \\ \lambda^* \end{bmatrix} = P\bar{z}$.

Recall that multiplication with a permutation matrix ($P$ and $P^T$) is done by arranging the elements in the vector. Matrix $B$ is 1 or 2 dimensional block diagonal therefore computing $\hat{z}$ is inexpensive. Cost of triangular substitutions with $L$ and $L^T$ depends on the sparsity of $L$. Normally, the significant cost for solving the system comes from the cost of performing factorization and triangular substitution, the latter of which depends on the sparsity of $L$. Observe that $B$ is a block diagonal matrix
$$B = \begin{bmatrix} B^{(1)} & 0 & \cdots & 0 \\ 0 & B^{(2)} & \cdots & 0 \\ \vdots & \vdots & \ddots & \vdots \\ 0 & 0 & \cdots & B^{(T)} \end{bmatrix},$$
where blocks $B^{(t)}$ are either 1-by-1 or 2-by-2 matrix and nonsingular. To perform symmetric indefinite factorization, let $K^{(t)}$ be the matrix that remains to be factorized in the $t$th iteration. The algorithm starts with $K^{(1)} = K$. For each iteration, we first identify a submatrix $B^{(t)}$ from elements of $K^{(t)}$ that are suitable to be used as a pivot block (There are many methods for selecting a suitable pivot $B^{(t)}$. Our method is described in Section 2.3). The submatrix $B^{(t)}$ is either a single diagonal element of $K^{(t)}$ ($\begin{bmatrix} k_{ll}^{(t)} \end{bmatrix}$) or a 2-by-2 block with two diagonal elements of $K^{(t)}$ $\left( \begin{bmatrix} k_{ll}^{(t)} & k_{lr}^{(t)} \\ k_{rl}^{(t)} & k_{rr}^{(t)} \end{bmatrix} \right)$.

Next, we find the permutation matrix $P^{(t)}$ satisfying
$$\left(P^{(t)}\right)^T K^{(t)} P^{(t)} = \begin{bmatrix} B^{(t)} & \left(C^{(t)}\right)^T \\ C^{(t)} & Z^{(t)} \end{bmatrix}. \quad (5)$$

The right-hand side of (5) can be factorized as
$$\left(P^{(t)}\right)^T K^{(t)} P^{(t)} = \begin{bmatrix} I & 0 \\ C^{(t)}\left(B^{(t)}\right)^{-1} & I \end{bmatrix} \cdot$$
$$\begin{bmatrix} B^{(t)} & 0 \\ 0 & Z^{(t)} - C^{(t)}\left(B^{(t)}\right)^{-1}\left(C^{(t)}\right)^T \end{bmatrix} \cdot$$
$$\begin{bmatrix} I & \left(B^{(t)}\right)^{-1}\left(C^{(t)}\right)^T \\ 0 & I \end{bmatrix}. \quad (6)$$

Let $L^{(t)} = C^{(t)}\left(B^{(t)}\right)^{-1}$ and $K^{(t+1)} = Z^{(t)} - C^{(t)}\left(B^{(t)}\right)^{-1}\left(C^{(t)}\right)^T$. The above can be rewritten as
$$\left(P^{(t)}\right)^T K^{(t)} P^{(t)} = \begin{bmatrix} I & 0 \\ L^{(t)} & I \end{bmatrix} \cdot \begin{bmatrix} B^{(t)} & 0 \\ 0 & K^{(t+1)} \end{bmatrix} \cdot$$
$$\begin{bmatrix} I & \left(L^{(t)}\right)^T \\ 0 & I \end{bmatrix}. \quad (7)$$

The same process can be repeated recursively on the matrix $K^{(t+1)}$. Note that the dimension of $K^{(t+1)}$ is less than the dimension of $K^{(t)}$ by either one or two depending on the dimension of $B^{(t)}$. Choosing pivot at each step should be inexpensive, lead to at most modest growth in the elements of the remaining matrix, and $L$ should not be too much denser than the original matrix. There are various method to identify pivot block $B^{(t)}$ for dense matrices. Bunch and Parlett searches the whole submatrix at each stage for the largest-magnitude diagonal $k_{qq}^{(t)}$ and the largest-magnitude off-diagonal $k_{rl}^{(t)}$. It identifies $k_{qq}^{(t)}$ as the 1-by-1 pivot block if the resulting growth rate is acceptable. Otherwise, it selects $\begin{bmatrix} k_{ll}^{(t)} & k_{lr}^{(t)} \\ k_{rl}^{(t)} & k_{rr}^{(t)} \end{bmatrix}$ as the 2-by-2 pivot block. This method requires $O(n^3)$ comparisons and yields a matrix $L$ whose maximum element is bounded by 2.781. Bunch-Kaufman pivoting strategy searches for the largest-magnitude off-diagonal elements of at most two columns for each iteration. It requires $O(n^2)$ comparisons but the elements in $L$ are unbounded. BBK combines the two above strategies and is widely used to select pivot blocks. By monitoring the size of the elements in $L$, BBK uses the Bunch-Kaufman strategy when it yields modest element growth. Otherwise, it repeatedly searches for an acceptable pivot [13]. BBK algorithm is shown in Algorithm 1 below. On average cases, the total cost of BBK is the same as Bunch-Kaufman, but in the worst case it can be the same as the cost of the Bunch-Parlett strategy.

## 2.3 A New Pivot Selection for Block Constraint Quadratic Programming

This section describes our proposed pivot selection method for the quadratic programs with block diagonal constraint matrices. The goals of our method are to maintain sparsity and stability in the factors. First, we identify candidate pivots that can maintain sparsity of $L$. Second, we select among these candidates to maintain stability of the factors. The last subsection describes the overall algorithm for factoring the KKT matrix.

**Algorithm 1** The BBK algorithm
Set $\alpha = (1 + \sqrt{17})/8$
Set $\gamma_1$ = maximum magnitude of any subdiagonal entry in column 1
**if** $|k_{11}| \geq \alpha\gamma_1$ **then**
  Use $k_{11}$ as a 1×1 pivot
**else**
  Set $l = 1; \gamma_l = \gamma_1$ ;
  **repeat**
    Set $r$ = row index of first (subdiagonal) entry of maximum magnitude in column $l$
    Set $\gamma_r$ = maximum magnitude of any off-diagonal entry in column $r$
    **If** $|k_{rr}| \geq \alpha\gamma_r$ **then**
      Use $k_{rr}$ as a 1×1 pivot
    **else if** $\gamma_l = \gamma_r$ **then**
      Use $\begin{bmatrix} k_{ll} & k_{lr} \\ k_{rl} & k_{rr} \end{bmatrix}$ as 2×2 pivot
    **else**
      Set $l = r; \gamma_l = \gamma_r$
    **end if**
  **until** A pivot is chosen
**end if**

### 2.3.1 Candidate Pivots Identification

Consider the structure of the KKT matrix of our quadratic program with a block diagonal constraint matrix. Elements of the KKT matrix $K$ can be classified into three types: the elements of the Hessian matrix $h_{ij}$, the nonzero elements of the constraint matrix $a_{ij}$, and the zero submatrices. For better readability, we use a 6-by-6 Hessian matrix with two blocks of constraints as an example in our explanation. The structure of a sample KKT matrix is as follows:

$K =$
$$\begin{bmatrix} h_{11} & h_{12} & h_{13} & h_{14} & h_{15} & h_{16} & a_{11} & a_{21} & 0 & 0 \\ h_{21} & h_{22} & h_{23} & h_{24} & h_{25} & h_{26} & a_{12} & a_{22} & 0 & 0 \\ h_{31} & h_{32} & h_{33} & h_{34} & h_{35} & h_{36} & a_{13} & a_{23} & 0 & 0 \\ h_{41} & h_{42} & h_{43} & h_{44} & h_{45} & h_{46} & 0 & 0 & a_{34} & a_{44} \\ h_{51} & h_{52} & h_{53} & h_{54} & h_{55} & h_{56} & 0 & 0 & a_{35} & a_{45} \\ h_{61} & h_{62} & h_{63} & h_{64} & h_{65} & h_{66} & 0 & 0 & a_{36} & a_{46} \\ a_{11} & a_{12} & a_{13} & 0 & 0 & 0 & 0 & 0 & 0 & 0 \\ a_{21} & a_{22} & a_{23} & 0 & 0 & 0 & 0 & 0 & 0 & 0 \\ 0 & 0 & 0 & a_{34} & a_{35} & a_{36} & 0 & 0 & 0 & 0 \\ 0 & 0 & 0 & a_{44} & a_{45} & a_{46} & 0 & 0 & 0 & 0 \end{bmatrix}.$$
(8)

Note that $h_{ij} = h_{ji}$ due to $H$ being symmetric. For our KKT matrix, there are three possible cases for pivot $B^{(t)}$. The first case is a 1-by-1 matrix selected from one of the nonzero diagonal elements in matrix $K^{(t)}$ (i.e., $h_{ll}^{(t)}$). The second case is a 2-by-2 matrix where both diagonal elements are nonzero (i.e., $\begin{bmatrix} h_{ll}^{(t)} & h_{lr}^{(t)} \\ h_{rl}^{(t)} & h_{rr}^{(t)} \end{bmatrix}$). In this case, the off-diagonal elements are the elements of the Hessian matrix $h_{ij}$, where $i \neq j$. The last possible case is a 2-by-2 matrix where one diagonal element is zero and the other three elements are nonzero (i.e., $\begin{bmatrix} h_{ll}^{(t)} & a_{rl}^{(t)} \\ a_{rl}^{(t)} & 0 \end{bmatrix}$). In other words, the off-diagonal elements are the nonzero elements of the constraint matrix $a_{ij}$. Selecting a pivot in any other ways besides the three mentioned above is not possible as they all lead to singular $B^{(t)}$. Each form of pivot $B^{(t)}$ directly affects the sparsity of the factor $L^{(t)}$ and also the sparsity and the stability of the remaining matrix $K^{(t+1)}$. Note that the sparsity of $K^{(t+1)}$ affects the sparsity of $L^{(t+1)}$, too (Recall that $L^{(t)} = C^{(t)}(B^{(t)})^{-1}$ and $K^{(t+1)} = Z^{(t)} - C^{(t)}(B^{(t)})^{-1}(C^{(t)})^T$). Now we consider the three cases of pivot in more details. For the first case, where $B^{(t)}$ is a 1-by-1 matrix, the number of zeros in $L^{(t)}$ is equal to the number of zeros in $C^{(t)}$ but many zeros in $K^{(t)}$ become nonzeros (fill-ins) in the remaining matrix $K^{(t+1)}$. Selecting pivot of this form generally cannot maintain the sparsity of the factors. For example, let $B^{(1)} = [h_{55}]$ be the pivot for the matrix in (8). After permutation, we have

$(P^{(1)})^T K^{(1)} P^{(1)} =$
$\begin{bmatrix} h_{55} & h_{52} & h_{53} & h_{54} & h_{51} & h_{56} & 0 & 0 & a_{35} & a_{45} \\ h_{25} & h_{22} & h_{23} & h_{24} & h_{21} & h_{26} & a_{12} & a_{22} & 0 & 0 \\ h_{35} & h_{32} & h_{33} & h_{34} & h_{31} & h_{36} & a_{13} & a_{23} & 0 & 0 \\ h_{45} & h_{42} & h_{43} & h_{44} & h_{41} & h_{46} & 0 & 0 & a_{34} & a_{44} \\ h_{15} & h_{12} & h_{13} & h_{14} & h_{11} & h_{16} & a_{11} & a_{21} & 0 & 0 \\ h_{65} & h_{62} & h_{63} & h_{64} & h_{61} & h_{66} & 0 & 0 & a_{36} & a_{46} \\ 0 & a_{12} & a_{13} & 0 & a_{11} & 0 & 0 & 0 & 0 & 0 \\ 0 & a_{22} & a_{23} & 0 & a_{21} & 0 & 0 & 0 & 0 & 0 \\ a_{35} & 0 & 0 & a_{34} & 0 & a_{36} & 0 & 0 & 0 & 0 \\ a_{45} & 0 & 0 & a_{44} & 0 & a_{46} & 0 & 0 & 0 & 0 \end{bmatrix}$,

$C^{(1)} = \begin{bmatrix} h_{25} \\ h_{35} \\ h_{45} \\ h_{15} \\ h_{65} \\ 0 \\ 0 \\ a_{35} \\ a_{45} \end{bmatrix}$, $L^{(1)} = \begin{bmatrix} \times \\ \times \\ \times \\ \times \\ \times \\ 0 \\ 0 \\ \times \\ \times \end{bmatrix}$, and

$K^{(2)} = \begin{bmatrix} \times & \times & \times & \times & \times & \times & \times & \bullet & \bullet \\ \times & \times & \times & \times & \times & \times & \times & \bullet & \bullet \\ \times & \times & \times & \times & \times & 0 & 0 & \times & \times \\ \times & \times & \times & \times & \times & \times & \times & \bullet & \bullet \\ \times & \times & \times & \times & \times & 0 & 0 & \times & \times \\ \times & \times & 0 & \times & 0 & 0 & 0 & \bullet & \bullet \\ \times & \times & 0 & \times & 0 & 0 & 0 & \bullet & \bullet \\ \bullet & \bullet & \times & \bullet & \times & \bullet & \bullet & \bullet & \bullet \\ \bullet & \bullet & \times & \bullet & \times & \bullet & \bullet & \bullet & \bullet \end{bmatrix}$.

Note that $\times$ denotes the nonzero elements in $L^{(1)}$ and $K^{(2)}$ in which the elements of matrix $C^{(1)}$ and $Z^{(1)}$ (in the same positions) are also nonzero and $\bullet$ denotes the fill-in elements compared to $C^{(1)}$ and $Z^{(1)}$, respectively. For the second case where the pivot is a 2-by-2 matrix with no zero elements, the number of nonzeros in $L^{(t)}$ may be equal to or greater than that of $C^{(t)}$. It also results in a large number of fill-ins in the remaining matrix $K^{(t+1)}$. For example, suppose $B^{(1)} = \begin{bmatrix} h_{55} & h_{56} \\ h_{65} & h_{66} \end{bmatrix}$. We have

$(P^{(1)})^T K^{(1)} P^{(1)} =$
$\begin{bmatrix} h_{55} & h_{56} & h_{53} & h_{54} & h_{51} & h_{52} & 0 & 0 & a_{35} & a_{45} \\ h_{65} & h_{66} & h_{63} & h_{64} & h_{61} & h_{62} & 0 & 0 & a_{36} & a_{46} \\ h_{35} & h_{36} & h_{33} & h_{34} & h_{31} & h_{32} & a_{13} & a_{23} & 0 & 0 \\ h_{45} & h_{46} & h_{43} & h_{44} & h_{41} & h_{42} & 0 & 0 & a_{34} & a_{44} \\ h_{15} & h_{16} & h_{13} & h_{14} & h_{11} & h_{12} & a_{11} & a_{21} & 0 & 0 \\ h_{25} & h_{26} & h_{23} & h_{24} & h_{21} & h_{22} & a_{12} & a_{22} & 0 & 0 \\ 0 & 0 & a_{13} & 0 & a_{11} & a_{12} & 0 & 0 & 0 & 0 \\ 0 & 0 & a_{23} & 0 & a_{21} & a_{22} & 0 & 0 & 0 & 0 \\ a_{35} & a_{36} & 0 & a_{34} & 0 & 0 & 0 & 0 & 0 & 0 \\ a_{45} & a_{46} & 0 & a_{44} & 0 & 0 & 0 & 0 & 0 & 0 \end{bmatrix}$,

$C^{(1)} = \begin{bmatrix} h_{35} & h_{36} \\ h_{45} & h_{46} \\ h_{15} & h_{16} \\ h_{25} & h_{26} \\ 0 & 0 \\ 0 & 0 \\ a_{35} & a_{36} \\ a_{45} & a_{46} \end{bmatrix}$, $L^{(1)} = \begin{bmatrix} \times & \times \\ \times & \times \\ \times & \times \\ \times & \times \\ 0 & 0 \\ 0 & 0 \\ \times & \times \\ \times & \times \end{bmatrix}$, and

$K^{(2)} = \begin{bmatrix} \times & \times & \times & \times & \times & \times & \bullet & \bullet \\ \times & \times & \times & \times & 0 & 0 & \times & \times \\ \times & \times & \times & \times & \times & \times & \bullet & \bullet \\ \times & \times & \times & \times & \times & \times & \bullet & \bullet \\ \times & 0 & \times & \times & 0 & 0 & 0 & 0 \\ \times & 0 & \times & \times & 0 & 0 & 0 & 0 \\ \bullet & \times & \bullet & \bullet & 0 & 0 & \bullet & \bullet \\ \bullet & \times & \bullet & \bullet & 0 & 0 & \bullet & \bullet \end{bmatrix}$.

Note that $K^{(t+1)}$ can be denser than in the above example for some other pivots such as $B^{(1)} = \begin{bmatrix} h_{33} & h_{36} \\ h_{63} & h_{66} \end{bmatrix}$. Lastly, consider the case where the pivot is a 2-by-2 matrix with one diagonal element being zero. In this case, the number of zeros in $L^{(t)}$ is equal to the number of zeros in $C^{(t)}$ and there is no fill-in in the remaining matrix $K^{(t+1)}$. For example, suppose $B^{(1)} = \begin{bmatrix} h_{55} & a_{35} \\ a_{35} & 0 \end{bmatrix}$. We have

$(P^{(1)})^T K^{(1)} P^{(1)} =$

$$\begin{bmatrix} h_{55} & a_{35} & h_{53} & h_{54} & h_{51} & h_{65} & 0 & 0 & h_{52} & a_{45} \\ a_{35} & 0 & 0 & a_{34} & 0 & a_{36} & 0 & 0 & 0 & 0 \\ h_{53} & 0 & h_{33} & h_{43} & h_{31} & h_{63} & a_{13} & a_{23} & h_{32} & 0 \\ h_{54} & a_{34} & h_{43} & h_{44} & h_{41} & h_{64} & 0 & 0 & h_{42} & a_{44} \\ h_{51} & 0 & h_{31} & h_{41} & h_{11} & h_{61} & a_{11} & a_{21} & h_{21} & 0 \\ h_{65} & a_{36} & h_{63} & h_{64} & h_{61} & h_{66} & 0 & 0 & h_{62} & a_{46} \\ 0 & 0 & a_{13} & 0 & a_{11} & 0 & 0 & 0 & a_{12} & 0 \\ 0 & 0 & a_{23} & 0 & a_{21} & 0 & 0 & 0 & a_{22} & 0 \\ h_{52} & 0 & h_{32} & h_{42} & h_{21} & h_{62} & a_{12} & a_{22} & h_{22} & 0 \\ a_{45} & 0 & 0 & a_{44} & 0 & a_{46} & 0 & 0 & 0 & 0 \end{bmatrix},$$

$$C^{(1)} = \begin{bmatrix} h_{53} & 0 \\ h_{54} & a_{34} \\ h_{51} & 0 \\ h_{65} & a_{36} \\ 0 & 0 \\ 0 & 0 \\ h_{52} & 0 \\ a_{45} & 0 \end{bmatrix}, L^{(1)} = \begin{bmatrix} 0 & \times \\ \times & \times \\ 0 & \times \\ \times & \times \\ 0 & 0 \\ 0 & 0 \\ 0 & \times \\ 0 & \times \end{bmatrix}, \text{ and }$$

$$K^{(2)} = \begin{bmatrix} \times & \times & \times & \times & \times & \times & \times & 0 \\ \times & \times & \times & \times & 0 & 0 & \times & \times \\ \times & \times & \times & \times & \times & \times & \times & 0 \\ \times & \times & \times & \times & 0 & 0 & \times & \times \\ \times & 0 & \times & 0 & 0 & 0 & \times & 0 \\ \times & 0 & \times & 0 & 0 & 0 & \times & 0 \\ \times & \times & \times & \times & \times & \times & \times & 0 \\ 0 & \times & 0 & \times & 0 & 0 & 0 & 0 \end{bmatrix}.$$

We see that the first and second types of pivots generate fill-in in remaining matrix $K^{(t+1)}$. The third type yields sparser $L^{(t)}$ than the other two and generally produces no fill-ins in the remaining matrix $K^{(t+1)}$. Therefore, our algorithm first identifies all candidate pivots that are of the form $\begin{bmatrix} h_{ii} & a_{ji} \\ a_{ji} & 0 \end{bmatrix}$. By choosing this form of pivots, the factor $L$ is as sparse as possible and there are no fill-ins in the remaining matrix $K^{(t+1)}$.

### 2.3.2 Pivot Selection

There are generally many pivot candidates of the form that we are interested in. We compare the condition numbers of these candidate pivots and then choose the one with the smallest condition number. Recall that the condition number of a 2-by-2 matrix is defined as
$$\text{cond}(B) = \|B\| \cdot \|B^{-1}\|. \qquad (9)$$

When the candidate pivot $B$ is of the form $\begin{bmatrix} b_{ii} & b_{ij} \\ b_{ij} & 0 \end{bmatrix}$, $B^{-1}$ becomes $\left(-\frac{1}{b_{ij}^2} \cdot \begin{bmatrix} 0 & -b_{ij} \\ -b_{ij} & b_{ii} \end{bmatrix}\right)$. Using infinity norm, we see that

$$\|B\|_\infty = \max\{|b_{ii}| + |b_{ij}|, |b_{ij}|\}$$
$$= |b_{ii}| + |b_{ij}|$$
$$\|B^{-1}\|_\infty = \left(\frac{1}{b_{ij}^2}\right) \max\{|b_{ij}|, |b_{ij}| + |b_{ii}|\}$$
$$= \frac{|b_{ij}| + |b_{ii}|}{b_{ij}^2}$$
$$\text{cond}_\infty(B) = \frac{(|b_{ii}| + |b_{ij}|) \cdot (|b_{ij}| + |b_{ii}|)}{b_{ij}^2}$$
$$= \left(1 + \frac{|b_{ii}|}{|b_{ij}|}\right)^2 \qquad (10)$$

Therefore, we need only to compare $|b_{ii}|/|b_{ij}|$ to find the pivot candidate with the minimum condition number. We do so and select the candidate with the smallest condition number as the pivot. Note that, when $a_{ij}$ is zero, the condition numbers of the candidates containing this $a_{ij}$ are infinity. In this case, such candidates are not chosen by our algorithm.

### 2.3.3 The Algorithm

This section gives the complete picture of our algorithm. First, we choose the pivot that maintains the sparsity of the factors. We select a 2-by-2 pivot matrix with one of the diagonal elements being zero as described previously. These pivots yield no fill-ins and we can choose this type of pivots for the first $m$ iterations, where $m$ is the number of constraints in the quadratic program. Afterward, $K^{(m+1)}$ is completely dense therefore we switch to use a general (non-sparse) symmetric indefinite factorization at this point.

Our method keeps track of the current and original positions of elements $a_{ij}$ (as the elements may change in the permutation step). These positions are used to efficiently produce the pivot candidates of the third form. We consider only the candidates with the off-diagonal entries from the same block $A_i$, where $i$ is chosen arbitrarily. (The pivot candidates are identified from the elements from each block by block.) Among them, we select the candidate with the smallest condition number. Note that, according to (10), we need to compute only the condition numbers of the candidates with the largest $|b_{ij}|$ for each column $j$ and selecting the one with the smallest condition number. If pivot $\begin{bmatrix} k_{ll}^{(t)} & k_{lr}^{(t)} \\ k_{rl}^{(t)} & 0 \end{bmatrix}$ is selected as the pivot, we remove row $r$ and column $l$ from the lists of available row and column. If row $r$ is the last row in a block, we remove the block containing row $r$ from the available block list, too. Then we continue to the pivot candidates from the next block. Note that our pivot selection method requires $O(\sum_{i=1}^{N} m_i^2 n_i)$ operations. Our method is shown in Algorithm 2 below.

## 3. RESULTS AND DISCUSSION

In this section, we compare the efficiency between the following two methods: (i) MA57 and (ii) our method. The experiment was performed in Matlab 2012a on problems with 500, 1000, and 1500 variables. For each problem size, we test with 10, 50, and 100 blocks in the constraint matrix, where each block is of equal size. The numbers of constraints are 40 and 80 percent of the number of variables. The test problems are randomly generated in the following way: Let $\widehat{H} \in \mathbb{R}^{n \times n}$. Each element of $\widehat{H}$, $c$, and $e$, and each nonzero element of $A$ is randomly generated between zero and one according to the uniform distribution. Then, let $H = \widehat{H}\widehat{H}^T / (\max_{i,j} \hat{h}_{i,j})$. For each problem, we experiment with 10 different instances. We compare the average numbers of nonzeros in factor $L$. The results of this experiment, which are shown in Table 1, show that

**Algorithm 2** Symmetric indefinite factorization for QP block constraint KKT matrix
Set $K$ = KKT matrix of QP with block constraints
Set $N$ = number of blocks in constraints, $n$ = number of all constraints
Set $n$ = number of variables, $s = n + m$ // size of matrix $K$
Set $L = s$-by-$s$ identity matrix, $B = s$-by-$s$ zero matrix
Set $aB = \{1, 2, ..., \}$ // list of available block
Set $aC = \{1, 2, ..., n\}$ // list of available column
Set $aR = \{n + 1, n + 2, ..., \}$ // list of available row
Set $P = [1\ 2\ ...\ s]$ // list of columns (1 to $n$) and rows ($n + 1$ to $s$) position in matrix
Set $sR = [sR_1, sR_2, ..., sR_N]$ // $sR_i$ is the first row of $A_i$
Set $eR = [eR_1, eR_2, ..., eR_N]$ // $eR_i$ is the last row of $A_i$
Set $sC = [sC_1, sC_2, ..., sC_N]$ // $sC_i$ is the first column of $A_i$
Set $eC = [eC_1, eC_2, ..., eC_N]$ // $eC_i$ is the last column of $A_i$
Set $p = 1$
**while** $p < m \times 2$ **do**
    Set $mincond = \infty$, $removeBl = 0$
    Randomly select $t$ from $aB$
    Set $avaiRowInBl = \{x : x \in aR; sR_t \leq x \leq eR_t\}$
    Set $avaiColInBl = \{x : x \in aC; sC_t \leq x \leq eC_t\}$
    Set $posRowInBl = \{x : x = P_i; i \in avaiRowInBl\}$
    Set $posColInBl = \{x : x = P_i; i \in avaiColInBl\}$
    Set $mincond = \min\left\{\left(\left|\frac{K_{jj}}{K_{ij}}\right|\right) : i \in posRowInBl, j \in posColInBl\right\}$
    Set $l$ = column of $mincond$
    Set $r$ = row of $mincond$
    **if** $|avaiRowInBl| = 1$ **then**
        Remove $t$ from $aB$
    **end if**
    Use $\begin{bmatrix} k_{ll} & k_{lr} \\ k_{lr} & k_{rr} \end{bmatrix}$ as the 2-by-2 pivot
    Remove elements $\{x : x = P_r \text{ or } x = P_l\}$ from $aR$ and $aC$
    Swap $P_l$ and $P_r$ to $P_p$ and $P_p + 1$, respectively
    $p = p + 2$
**end while**
Factorize the remaining matrix with a general (non-sparse) symmetric indefinite factorization method our method yields sparser *L* when computing the symmetric indefinite factorization than the MA57 algorithm. After the factorization, we use the factors from the two methods to compute the solution of the quadratic program following Steps (i) - (iv) in Section 2.2. Table 1 also shows the solving time and the accuracy of our algorithm. The results show that using the factors *L*, *B*, and *P* from our method reduces the time needed to solve the KKT system compared to using the factors from MA57. Our method also yields accurate solutions with small residuals.

Normally, the Hessian matrix may not be dense. We therefore also experiment on the problems with sparse Hessian matrices having 30, 50, and 70 percent of their entries being nonzeros. We test with the constraint matrices having 10, 50, and 100 blocks, where each block is of equal size. The numbers of constraints are 40 and 80 percent of the number of variables. For each problem, we experiment with 10 different instances. We compare the average numbers of nonzeros in the factor *L*. The results of this experiment are shown in Table 2. We see that even when the Hessian matrix is sparse, our method still maintains more sparsity in *L* than MA57 can.

Finally, we compare both methods on problems where each blocks in the constraint matrices are of different sizes. The results are shown in Table 3. For these problems, our method also produces sparser factors and requires less solving time than MA57.

## 4. CONCLUSIONS

In this article, we propose an algorithm for solving block constraint quadratic programs. Our method is a direct method using symmetric indefinite factorization. This method exploits the known structure of the quadratic problem to efficiently compute the factors that are stable and retain the sparsity of the problem. The results of our experiments show that our method is better at maintaining sparsity of the factors than the MA57 algorithm. Consequently, using the factors from our method to solve the KKT system is faster than using the factors from MA57 while yielding the solution that is as accurate. Finally, we note that the steps in our pivot selection algorithm are easily parallelized and therefore can be made more efficient with parallel computing, too.

**Acknowledgements**

This work is partially supported by Center of Excellence in Intelligent Informatics, Speech and Language Technology, and Service Innovation (CILS), Thammasat University and National Research Universities, Thailand.

**Table 1.** Average numbers of nonzeros in factor $L$, average solving time, and average residual of MA57 and our algorithm for problems with 500, 1000, and 1500 variables and constraint matrices with equal-sized blocks.

| $n$ | $N \times (n_i \times m_i)$ | Ave. num. of nonzeros in L | | Ave. solving time (ms) | | Ave. residual ($\times 10^{-10}$) | |
|---|---|---|---|---|---|---|---|
| | | MA57 | Our Method | MA57 | Our Method | MA57 | Our Method |
| 500 | 10×(50×10) | 162210.0 | 130250.0 | 3.81240 | 2.39233 | 0.046 | 0.028 |
| 500 | 10×(50×40) | 333090.0 | 145250.0 | 6.58479 | 4.37043 | 0.190 | 0.033 |
| 500 | 20×(25×5) | 158829.4 | 127750.0 | 3.43303 | 2.32784 | 0.042 | 0.032 |
| 500 | 20×(25×20) | 319564.8 | 135250.0 | 6.69930 | 4.31164 | 0.208 | 0.033 |
| 500 | 50×(10×2) | 156691.3 | 126250.0 | 3.65311 | 2.32088 | 0.037 | 0.030 |
| 500 | 50×(10×8) | 311022.7 | 129250.0 | 6.40970 | 4.35027 | 0.302 | 0.032 |
| 1000 | 10×(100×20) | 648420.0 | 520500.0 | 15.14728 | 7.80979 | 0.169 | 0.113 |
| 1000 | 10×(100×80) | 1332180.0 | 580500.0 | 32.11714 | 15.89429 | 0.620 | 0.138 |
| 1000 | 20×(50×10) | 634909.3 | 510500.0 | 15.63845 | 7.83631 | 0.146 | 0.113 |
| 1000 | 20×(50×40) | 1278127.7 | 540500.0 | 31.37790 | 15.77975 | 0.560 | 0.132 |
| 1000 | 50×(20×4) | 626351.3 | 504500.0 | 15.16187 | 7.79729 | 0.128 | 0.109 |
| 1000 | 50×(20×16) | 1243965.5 | 516500.0 | 30.59341 | 16.00238 | 0.478 | 0.131 |
| 1500 | 10×(150×30) | 1458630.0 | 1170750.0 | 37.22620 | 16.65317 | 0.336 | 0.260 |
| 1500 | 10×(150×120) | 2997261.7 | 1305750.0 | 83.14423 | 35.26495 | 1.104 | 0.295 |
| 1500 | 20×(75×15) | 1428238.8 | 1148250.0 | 36.45427 | 16.53619 | 0.300 | 0.227 |
| 1500 | 20×(75×60) | 2875679.8 | 1215750.0 | 80.20526 | 34.66840 | 1.019 | 0.276 |
| 1500 | 50×(30×6) | 1408997.6 | 1134750.0 | 35.57901 | 16.43024 | 0.286 | 0.231 |
| 1500 | 50×(30×24) | 2798788.2 | 1161750.0 | 75.85210 | 34.51535 | 1.105 | 0.268 |

Column $n$ is the number of variables. Column $N(n_i \times m_i)$ indicates the dimensions of each diagonal blocks in the constraint matrix. Columns Ave. num. of nonzeros in $L$ show the number of nonzeros in $L$ from MA57 algorithm and our method, respectively. Columns Ave. solving time show the solving time of the two methods. Columns Ave. residual represent the residuals of the results of both methods. The residual is $\|Kx - v\|_2$, where $K$ is our KKT matrix, $x$ is the computed solution, and $v$ is the vector $\begin{bmatrix} g \\ h \end{bmatrix}$ in (3).

**Table 2.** Average numbers of nonzeros in $L$ of MA57 algorithm and our algorithm for problems with 1000 variables with 30, 50, and 70% of nonzeros in Hessian matrix and constraint matrices with equal-sized blocks.

| | | Ave. num. of nonzeros in $L$ | | | | | |
|---|---|---|---|---|---|---|---|
| | | 30% nonzeros in $H$ | | 50% nonzeros in $H$ | | 70% nonzeros in $H$ | |
| $n$ | $N \times (n_i \times m_i)$ | MA57 | Our Method | MA57 | Our Method | MA57 | Our Method |
| 1000 | 10×(100×20) | 631400 | 507800 | 643900 | 514800 | 640900 | 518000 |
| 1000 | 10×(100×80) | 1509800 | 567900 | 1524700 | 575000 | 1441600 | 578100 |
| 1000 | 20×(50×10) | 612900 | 486700 | 622400 | 499900 | 632300 | 505800 |
| 1000 | 20×(50×40) | 1295700 | 516300 | 1296600 | 530100 | 1273900 | 536000 |
| 1000 | 50×(20×4) | 593800 | 455200 | 614300 | 479900 | 620900 | 493500 |
| 1000 | 50×(20×16) | 1320100 | 458100 | 1301300 | 491300 | 1279600 | 505600 |

Column $n$ is the number of variables. $N(n_i \times m_i)$ indicates the dimensions of each diagonal blocks in the constraint matrix. Columns Ave. num. of nonzeros in $L$ show the number of nonzeros in $L$ from MA57 algorithm and our method, respectively.

**Table 3.** Average numbers of nonzeros in factor $L$, average solving time, and average residual of MA57 and our algorithm for problems with 500, 1000, and 1500 variables constraint matrices with unequal-sized blocks.

| | | Ave. num. of nonzeros in $L$ | | Ave. solving time (ms) | | Ave. residual ($\times 10^{-10}$) | |
|---|---|---|---|---|---|---|---|
| $n$ | $N \times (n_i \times m_i)$ | MA57 | Our Method | MA57 | Our Method | MA57 | Our Method |
| 500 | 5× (50×10,75×15, 100×20,125×25, 150×30) | 160937.0 | 136500.0 | 3.54694 | 2.33908 | 0.053 | 0.033 |
| 500 | 5× (50×40,75×60, 100×80,125×100, 150×120) | 328058.0 | 170250.0 | 6.59455 | 4.74211 | 0.221 | 0.035 |
| 1000 | 5×(100×20,150×30, 200×40,250×50, 300×60) | 643154.0 | 545500.0 | 14.98070 | 7.85096 | 0.186 | 0.142 |
| 1000 | 5×(100×80,150×120, 200×160,250×200, 300×240) | 1311324.0 | 680500.0 | 31.48018 | 16.05201 | 0.824 | 0.155 |
| 1500 | 5×(150×30,225×45, 300×60,375×75, 450×90) | 1446635.9 | 1227000.0 | 37.44730 | 16.52353 | 0.468 | 0.420 |
| 1500 | 5×(150×120,225×180, 300×240,375×300, 450×360) | 2949798.0 | 1530750.0 | 81.43578 | 35.06723 | 1.283 | 0.610 |

Column $n$ is the number of variables. Column $N \times (n_i \times m_i)$ indicates the dimensions of each diagonal blocks in the constraint matrix. Columns Ave. num. of nonzeros in $L$ show the number of nonzeros in $L$ from MA57 algorithm and our method, respectively. Columns Ave. solving time show the solving time of the two methods. Columns Ave. residual represent the residuals of the results of both methods. The residual is $\|Kx - v\|_2$, where $K$ is our KKT matrix, $x$ is the computed solution, and $v$ is the vector $\begin{bmatrix} g \\ h \end{bmatrix}$ in (3).